# On domains of elliptic operators with distributional coefficients

by Immanuel Zachhuber[a]

[a]. Department of Mathematics,
Freie Universität Berlin, Germany
Email: immanuel.zachhuber@fu-berlin.de

**Abstract**

We show how one can use recently gained insights from the study of singular SPDEs, more particularly the study of singular operators via the theory of Paracontrolled Distributions, to construct domains for (singular) elliptic operators. Formally we consider

$$A(u)\text{“=”}(1-\Delta)u+\nabla V\cdot\nabla u+\xi u+\mathrm{div}\,(\rho u),$$

where $V\in\mathcal{C}^{\delta}$, $\xi\in\mathcal{C}^{-2+\delta}$, $\rho\in\mathcal{C}^{-1+\delta}$, $\mathrm{div}\,\rho=0$ and which satisfy a structural assumption that is notably satisfied when $\xi$ is a *sub-critical noise*, see [19]. We also show that under this assumption, one can construct a continuous change of variables $\Theta$ which satisfies

$$A\Theta-(1-\Delta)\in\mathcal{L}(H^{2-\delta''};H^{\delta'})$$

which allows us to define $A$ rigorously and parametrise a domain. Moreover, for suitably regularised operators

$$A_{\varepsilon}(u):=(1-\Delta)u+\nabla V_{\varepsilon}\cdot\nabla u+(\xi_{\varepsilon}+c_{\varepsilon})\cdot u+\mathrm{div}\,(\rho_{\varepsilon}\cdot u),$$

we show that for a strongly converging regularised change of variables $\Theta_{\varepsilon}\to\Theta$ we have

$$A_{\varepsilon}\Theta_{\varepsilon}\to A\Theta\text{ in }\mathcal{L}(H^2;L^2)$$

which in particular implies norm resolvent convergence to a limiting closed operator.

Finally, we give a class of examples and show how to apply these results to prove strong analytical local well-posedness for a singular Schrödinger equation formally given by

$$i\partial_t u+(1-\Delta)u+\nabla V\cdot\nabla u+\xi\cdot u=-|u|^2u$$

for singular $V,\xi$ and that its solution is the limit of the solution of the classical solutions of

$$i\partial_t u_{\varepsilon}+(1-\Delta)u_{\varepsilon}+\nabla V_{\varepsilon}\cdot\nabla u_{\varepsilon}+(\xi+c_{\varepsilon})\cdot u_{\varepsilon}=-|u_{\varepsilon}|^2u_{\varepsilon}$$

with suitably chosen initial data.

## 1 Introduction

In this paper, we make sense of elliptic operators with distributional coefficients, formally

$$Au\text{“=”}(1-\Delta)u+\nabla V\cdot\nabla u+\xi u+\mathrm{div}\,(\rho u),\tag{1}$$

where $V\in\mathcal{C}^{\delta}$, $\xi\in\mathcal{C}^{-2+\delta}$, $\rho\in\mathcal{C}^{-1+\delta}$, $\mathrm{div}\,\rho=0$ and $\delta>0$ which also satisfy a structural assumption, see Assumption 1 and Appendix A for the definition of Hölder-Besov spaces $\mathcal{C}^{\alpha}$. Notable examples include when $\xi$ is a *sub-critical noise* in the language of [19] which includes coloured Gaussian noises; our result extends results like [11] and [2] in the case of the 2-/3- dimensional *Anderson Hamiltonian*, see also Example 3 for some details for the 2d Anderson Hamiltonian which motivates in part Assumption 1.

This is in principle a rather elementary problem and one sees immediately that as soon as either $\xi, V$ or $\rho$ is a genuine distribution, $A$ maps generic smooth functions into a space that is strictly worse than $L^2$, so any domain will not contain a core of smooth functions. In some specific cases, explicit domains have been constructed e.g. delta potentials, see [1]. Let us mention that in $d \geqslant 2$ the delta potential is too singular to be treated by the methods of this paper.

We employ tools from singular SPDEs, namely Paracontrolled Distributions [10], but surprisingly we obtain some results which hold in the entire subcritical regime without the need to employ the heavy machinery of Hairer's Regularity Structures [13] or the high-order Paracontrolled Calculus of Bailleul-Bernicot [4]. Rather, all we need is an exponential transform inspired by [14],[8] and [11] as well as elementary ideas from "first-order" paracontrolled calculus [10]. We should mention, however, that [19] (which provides the most interesting class of subcritical examples for us) used the BPHZ machinery for Regularity Structures to construct the stochastic terms and in the current article we always assume that the stochastic objects have already been constructed including renormalisation. In other words, our result is purely analytical in the language of singular SPDEs.

The first results on constructing explicit domains for this kind of singular operator(where $V = \rho = 0$ and $\xi$ being 2d periodic white noise and $V = 0 = \rho$) were in [2] using Paracontrolled Distributions, where also self-adjointness was proved. This was further developed in [11],[22] also using ideas from Paracontrolled Distributions and in [18] using Regularity Structures. See also [22] for a construction on manifolds, where also *Weyl asymptotics* were proven. We expect that with small modifications one can adapt that proof to our more general setting.

In [11], it was also observed that using only an exponential ansatz inspired by [14] and [8], see Section 2.1, one can construct a simpler form domain of $-\Delta + \xi$ for $\xi$ being $2/3-$dimensional white noise, see also [23] for a concise note on this. This was generalised and systematised to *subcritical noises* in [19] (with $V = \rho = 0$), where the self-adjoint operators were constructed as Friedrichs extensions.

In the case $\xi = 0$, the operator $A$ (1) is the formal generator of an SDE

$$dX_t = b(X_t)dt + dB_t$$

with distributional drift $b = \nabla V + \rho$, see [6],[17] and [9] amongst others. In particular, [6] includes a discussion of the Helmholtz decomposition of distributional divergence-free vector fields, which in part motivates the form of the gradient terms in our ansatz for $A$ in (1). They are also able to treat certain supercritical drifts.

Let us remark that the tools used in this article are comparatively simple– mainly integration by parts and a paraproduct-type decomposition– and while we formulate all results in the periodic setting, it is straightforward to extend them to smooth domains with suitable boundary conditions, compact manifolds, or $\mathbb{R}^d$ assuming some behaviour at $\infty$.

Another thing we prove in the current paper is that the change of variables (denoted by $\Theta$ in Theorem 2) which parametrises the domain of our singular operator $A$, actually even transforms $A$ into a lower-order perturbation of the Laplacian, i.e.

$$A\Theta - (1 - \Delta) \in \mathcal{L}(H^{2-\delta''}; H^{\delta'}),$$

for $\delta'', \delta' \ll 1$which means that intuitively, after a change of variables, the operator can be dealt with using classical tools. Note that the more natural transform for evolution equations is actually a conjugation with a change of variables, e.g. in [27, 24, 7] (in the case $V = \rho = 0$), one uses that there is a change of variables $\Theta$ s.t.

$$\Theta^{-1} A\Theta - (1 - \Delta) \in \mathcal{L}(H^{2-\delta''}; H^{\delta'})$$

which is then used to transform the heat/Schrödinger group via

$$\Theta^{-1}e^{tA}\Theta = e^{t\Theta^{-1}A\Theta} \quad \text{and } \Theta^{-1}e^{itA}\Theta = e^{it\Theta^{-1}A\Theta}$$

which allows to use that the conjugated group is a perturbation of the one with the Laplacian and transfer some bounds to the un-conjugated group using the invertibility of $\Theta$. The transfrom we construct in the current article does not immediately have that property, however one can likely tweak the construction similarly to [7].

However, using a simplified argument in the case $\rho = 0$ (which we call the *self-adjoint case*), we can still get a transform that allows us to prove local well-posedness of of strong solutions to Schrödinger, heat or wave equations with distributional $\xi, V$ see Theorem 8 and Section 6.

Below we state the set of assumptions we need to make in order to get the analysis to work and prove Theorem 2.

**Assumption 1. (Approximate solution of elliptic KPZ-type equation)** *Let $\delta > \delta' > 0$,*

*$\varepsilon \in (0,1)$ and we write $c_\varepsilon$ for some family of (possibly diverging in the regularisation parameter $\varepsilon$) constants.*

*We assume the following to hold;*

- *The smooth objects $\xi_\varepsilon, V_\varepsilon$, $\rho_\varepsilon$ converge in the following sense*

$$\begin{aligned} \xi_\varepsilon \to \xi \textit{ in } \mathcal{C}^{-2+\delta}(\mathbb{T}^d), \qquad & V_\varepsilon \to V \textit{ in } \mathcal{C}^{\delta}(\mathbb{T}^d) \\ \operatorname{div} \rho_\varepsilon = \operatorname{div} \rho = 0 \quad \textit{and} \quad & \rho_\varepsilon \to \rho \textit{ in } \mathcal{C}^{-1+\delta}(\mathbb{T}^d; \mathbb{R}^d) \end{aligned}$$

- *There exist approximate solutions $(W_\varepsilon, Z_\varepsilon)$ to an elliptic KPZ-type equation s.t.*

$$(1-\Delta)W_\varepsilon - |\nabla W_\varepsilon|^2 + \nabla W_\varepsilon \cdot \nabla V_\varepsilon + \xi_\varepsilon + c_\varepsilon =: Z_\varepsilon \tag{2}$$

$$Z_\varepsilon \to Z \quad \textit{in} \quad \mathcal{C}^{-1+\delta'}, \tag{3}$$

$$W_\varepsilon \to W \quad \textit{in} \quad \mathcal{C}^{\delta}. \tag{4}$$

- *The following singular objects converge as follows*

$$e^{V_\varepsilon + 2W_\varepsilon} Z_\varepsilon =: \tilde{Z}_\varepsilon \tag{5}$$

$$\tilde{Z}_\varepsilon \to \tilde{Z} = e^{V+2W} Z \textit{ in } \mathcal{C}^{-1+\delta'}(\mathbb{T}^d) \tag{6}$$

$$\rho_\varepsilon e^{V_\varepsilon + 2W_\varepsilon} \to \rho e^{V+2W} \textit{ in } \mathcal{C}^{-1+\delta'}(\mathbb{T}^d; \mathbb{R}^d) \tag{7}$$

*and*

$$e^{V_\varepsilon + 2W_\varepsilon} \nabla V_\varepsilon \cdot \rho_\varepsilon,\ e^{V_\varepsilon + 2W_\varepsilon} \nabla W_\varepsilon \cdot \rho_\varepsilon \to e^{V+2W} \nabla V \cdot \rho,\ e^{V+2W} \nabla W \cdot \rho \textit{ in } \mathcal{C}^{-1+\delta'}(\mathbb{T}^d). \tag{8}$$

The main theorem is given rigorously below.

**Theorem 2.** *Let $\delta > \delta' > 0$, $c_\varepsilon$ and $\xi_\varepsilon, V_\varepsilon$ etc as in Assumption 1 and define the operator*

$$A_\varepsilon(u) := (1-\Delta)u + \nabla V_\varepsilon \cdot \nabla u + (\xi_\varepsilon + c_\varepsilon) \cdot u + \operatorname{div}(\rho_\varepsilon u). \tag{9}$$

*Then there exists an invertible change of variables $\Theta_\varepsilon \in \mathcal{L}(H^2)$ which converges strongly in $\mathcal{L}(H^2; L^2)$ to an invertible limit $\Theta \in \mathcal{L}(H^2; \Theta(H^2))$ with $\Theta(H^2) \subset H^\delta$, dense in $L^2$, and we have that there exists a closed, densely defined operator $(A, \mathcal{D}(A))$ s.t.*

$$A_\varepsilon \Theta_\varepsilon - (1-\Delta) \to A\Theta - (1-\Delta) \quad \textit{in } \mathcal{L}(H^{2-\delta''}; H^{\delta'})$$

*for $\delta', \delta'' \ll 1$ which in turn implies that $A_\varepsilon \to A$ in the sense of norm resolvent convergence. Moreover, we may choose $\mathcal{D}(A) := \Theta(H^2)$ and we have the norm equivalence*

$$\|A\Theta u\|_{L^2} + \|\Theta u\|_{L^2} \approx \|u\|_{H^2}.$$

We give an example to motivate the appearance of the (possibly divergent) constants $c_\varepsilon$ in (2) and (9). Then we give some remarks on the restrictions put on the potentials $\xi, V,$ and $\rho$.

**Example 3. (2d Anderson Hamiltonian)** This is the "simplest" singular setting that requires *renormalisation,* see [2] for the first construction. We take $V = \rho = 0$ and $\xi$ is 2d periodic spatial white noise which satisfies formally $\mathbb{E}(\xi(x)\xi(y)) = \delta(x-y)$ and $\xi \in \mathcal{C}^{-1-\kappa}$ a.s. for $\kappa > 0$. This means it is just about not covered by the generic result from Section 5.1.

If $\xi_\varepsilon$ is a regularisation of the noise, the choice $W_\varepsilon = (1-\Delta)^{-1}\xi_\varepsilon$ leads us to try to make sense of the object $|\nabla W_\varepsilon|^2$ coming from the exponential ansatz which appears in (2). Since $W_\varepsilon \to W$ in $\mathcal{C}^{1-\kappa}$ which is (4), there is no reason to think that $|\nabla W_\varepsilon|^2$ converges since one can not apply nonlinear functions to distributions. Indeed, this sequence can be shown to have no well-defined limit but instead

$$: |\nabla W_\varepsilon| :^2 = |\nabla W_\varepsilon|^2 - \mathbb{E}(|\nabla W_\varepsilon|^2) = |\nabla W_\varepsilon|^2 - c_\varepsilon^{(2)}$$

does have a well-defined limit $: |\nabla W|^2 : \in \mathcal{C}^{-2\kappa}$, where $: \ldots :$ denotes *Wick ordering* and the expectation $\mathbb{E}(|\nabla W_\varepsilon|^2) = c_\varepsilon^{(2)} \sim \log(1/\varepsilon)$ diverges. Thus (2) reads

$$(1-\Delta) W_\varepsilon - |\nabla W_\varepsilon|^2 + \xi_\varepsilon + c_\varepsilon = - : |\nabla W_\varepsilon| :^2$$

and denoting by $Z_\varepsilon = - : |\nabla W_\varepsilon| :^2$ we have (3) and by our notation, $(W_\varepsilon, Z_\varepsilon)$ is an approximate solution to the KPZ-type equation (2).

Since the regularity of $W$–and thus $e^W$– is $C^{1-\kappa}$, the product $e^{W_\varepsilon} Z_\varepsilon$ converges to the limit $e^W Z$ in $\mathcal{C}^{-2\kappa}$ without additional input, i.e. (6). As soon as $\xi \in \mathcal{C}^{-\sigma}$ $\sigma > \frac{3}{2}$ (notably when $\xi$ is *3d periodic spatial white noise*), the $e^{W_\varepsilon} Z_\varepsilon$ does *not* automatically converge to a limit, but one needs to verify that additionally, see Lemma 2.40 in [11].

**Example 4. (Brox diffusion)** Brox [5] studied the SDE

$$d X_t = -\frac{1}{2} \dot{B}(X_t) \mathrm{d}t + d W_t$$

for a two-sided Brownian motion $B$ through its generator given by

$$L w = \frac{1}{2} e^B \partial_x (e^{-B} \partial_x w) = \frac{1}{2} \partial_{xx} w - \frac{1}{2} \partial_x B \partial_x w$$

so this operator is an example of the kind considered in our paper. Since $V = B \in \mathcal{C}^{\frac{1}{2}-\delta}$, $\xi = \rho = W = 0$ Assumption 1 is satisfied. Further, the use of the exponential is similar to the ansatz discussed in Section 2.1. See also [21] for a recent, related construction and [16] for a generalisation to Levy noises.

**Remark 5.** The assumptions on the coefficients $\xi, V$ are rather generic, assuming some subcritical negative regularity and the structural condition that the solution $W$ of the elliptic KPZ-type equation (2) exists. In the case $V = 0 = \rho$ this is what is constructed in [19], see also Section 5.2.

The assumptions on $\rho$ are a bit less natural: (7) is what one would expect by power counting, but (8) says that there should be some cancellation in the scalar product $\nabla(e^{V+2W})\cdot\rho$. In particular, a very simple special case of Theorem 2 is the case $\xi=V=W=0$ and $\rho\in\mathcal{C}^{-1+\delta}$ div $\rho=0$. We believe that if one uses the general machinery of high-order paracontrolled calculus [4], one can in fact replace (8) by a less restrictive subcritical assumption on $\rho$ but we decide not to pursue that here since it would make the paper considerably more technical.

In general, it seems likely that an analogous (perhaps more general) result can be shown using higher order paracontrolled calculus [4] or Regularity Structures on Besov spaces [15] but we decided to present this very light-weight approach where only integration by parts, chain rules and simple paraproduct decompositions appear.

We give two more results which are valid in the case $\rho\equiv 0$ which we call the "self-adjoint" case, a name which will be justified below.

**Theorem 6. (Self-adjointness)** *Consider $V,\xi,W$ as in Assumption 1 and $\rho=0$, then the resulting operator $(A,\mathcal{D}(A))$ is self-adjoint w.r.t. the weighted $L^2$ scalar product*

$$\langle u,v\rangle_V := \left(e^{\frac{V}{2}}u, e^{\frac{V}{2}}v\right)_{L^2} \tag{10}$$

*and has form domain $\mathcal{D}_V(\sqrt{A})=e^W H^1(\mathbb{T}^d)$ and domain $\mathcal{D}(A)=\Theta H^2(\mathbb{T}^d)$, where $\Theta$ is the map from Theorem 2.*

**Remark 7.** Clearly the norms $\|\cdot\|_{L^2}$ and $\|\cdot\|_{L^2_V} := \left\|e^{\frac{V}{2}}\cdot\right\|_{L^2}$ are equivalent since $e^{\pm V}\in L^\infty$ so for the domain $\mathcal{D}(A)$ it does not really matter whether we consider $A$ as a closed operator

$$A:\mathcal{D}_V(A)\subset L^2_V\to L^2_V \quad\text{or}\quad A:\mathcal{D}(A)\subset L^2\to L^2.$$

The difference, of course, is whether the operators are self-adjoint (or even symmetric) w.r.t. the scalar product and so the form domain

$$\mathcal{D}_V(\sqrt{A}) := \{f\in L^2_V: |\langle f, Af\rangle_V|<\infty\}$$

is in general different from

$$\mathcal{D}(\sqrt{A}) := \{f\in L^2: |(f, Af)_{L^2}|<\infty\}.$$

Using the self-adjointness together with knowledge of the domain of $A$ allows us to make sense of and prove well-posedness for a nonlinear Schrödinger equations (NLS) formally given by

$$i\partial_t u + (1-\Delta)u + \nabla V\nabla u + \xi u = -|u|^2 u \text{ on } \mathbb{T}^d \tag{11}$$

for $d\leqslant 3$ for initial data $u(0)=u_0$ in $\mathcal{D}(A)$ inspired by the approach from [11].

**Theorem 8. (Local well-posedness of stochastic NLS)** *Consider $V,\xi,W$ as in Assumption 1 and $\rho=0$ as well as $\Theta$ from Theorem 2 for dimension $d\leqslant 3$. Then, for every $u_0\in\mathcal{D}(A)$ there exists a time $T>0$ depending on $u_0$ and the other parameters s.t. there exists a unique mild solution to (11) i.e. which solves*

$$u(t)=e^{itA}u_0 + i\int_0^t e^{i(t-s)A}u|u|^2(s)\mathrm{d}s.$$

*Further* $u \in C([0,T];\mathcal{D}(A)) \cap C^1([0,T];L^2)$ *is unique, it depends continuously on* $u_0$ *and one has that it is the limit of the unique solution to the unique classical solution* $u_\varepsilon$ *to*

$$i\partial_t u_\varepsilon + (1-\Delta)u_\varepsilon + \nabla V_\varepsilon \nabla u_\varepsilon + \xi_\varepsilon u_\varepsilon + c_\varepsilon u_\varepsilon = -|u_\varepsilon|^2 u_\varepsilon \quad \textit{on } \mathbb{T}^d$$

*with initial data* $u_\varepsilon(0) = A_\varepsilon^{-1} A_\varepsilon u_0$ *in the following sense:*

$$\begin{cases} u_\varepsilon \to u & \textit{in } C_T L^2 \\ A_\varepsilon u_\varepsilon \to A u & \textit{in } C_T L^2 \\ \partial_t u_\varepsilon \to \partial_t u & \textit{in } C_T L^2 \end{cases}. \tag{12}$$

We give a few remarks on this theorem about limitations and comparisons to other known results.

**Remark 9.**

i. Since we only consider local-in-time solutions, it is not really important which sign the nonlinearity has (i.e. focussing/defocussing).

ii. From the proof we see that the $\mathcal{D}(A)$ norm controls the $L^\infty$ norm and hence we could consider a much larger class of nonlinearities $f(u)$ than the cubic polynomial one considered here. The nice feature of this type of nonlinearity is that one has a positive energy which controls the $e^W H^1$ norm of the solution so in particular controls $L_T^\infty L^p$ norms for

$$\begin{cases} p \in [1,\infty) & d \leqslant 2 \\ p \in \left[1, \frac{2d}{d-2}\right] & d \geqslant 3 \end{cases}.$$

We do not exploit this fact here but solutions to equations like (11) (with $V = 0$) with initial data in the form domain were studied in [11, 7, 25, 27].

iii. Relatedly, the restriction $d \leqslant 3$ comes from the observation that the embedding $H^2 \hookrightarrow L^\infty$ is only true in those dimensions. If one uses deeper results from dispersive PDE such as Strichartz estimates [24, 7], one can avoid the necessity for $L^\infty$ bounds so in principle one could obtain well-posedness results in higher dimensions with more work. If one is interested in $\xi$ being spatial white noise then $d = 4$ is the critical dimension and so it is anyway more natural to consider the equation in $d \leqslant 3$.

iv. If one does not choose this exact type of "prepared" initial data $u_\varepsilon(0) = A_\varepsilon^{-1} A_\varepsilon u_0$ but simply some other convergent $u_\varepsilon(0) \to u_0$, one still obtains that the solutions converge strongly in $C_T L^2$ but one does not expect the other two properties in (12) to hold.

**Notation:** We always work on the $d-$dimensional torus $\mathbb{T}^d = (\mathbb{R}/\mathbb{Z})^d$, writing $L^p = L^p(\mathbb{T}^d)$ to mean the standard $p-$integrable Lebesgue functions, $p \in [1,\infty]$, with norm $\|f\|_{L^p} = (\int_{\mathbb{T}^d} |f|^p dx)^{\frac{1}{p}}$ and $H^s = H^s(\mathbb{T}^d)$, $s \in \mathbb{R}$, for the $L^2$ based Sobolev spaces with norm $\|u\|_{H^s} = \left\|(1-\Delta)^{\frac{s}{2}} u\right\|_{L^2}$ and we moreover write $B_{p,q}^s$, $s \in \mathbb{R}, p, q \in [1,\infty]$ to mean Besov spaces, whose definition is recalled in Appendix A, which include the special case of *Hölder-Besov* spaces $\mathcal{C}^s = B_{\infty,\infty}^s$ which agree with the classical Hölder spaces $C^s$ if $s \in (0,1)$.

Further, for Banach spaces $X, Y$, we write

$$\mathcal{L}(X;Y) := \{B: X \to Y \text{ linear and bounded}\}$$

and $\mathcal{L}(X) := \mathcal{L}(X;X)$. Moreover, for an injective map $B \in \mathcal{L}(X;Y)$, we have

$$B \in \mathcal{L}(X, B(X)) \text{ invertible, where } \|f\|_{B(X)} := \|B^{-1} f\|_X$$

turns $B(X) \subset Y$ into a Banach space.

We also use the notations $\mathcal{F}, \mathcal{F}^{-1}$ to mean the Fourier and inverse Fourier transform on $\mathbb{T}^d$ respectively. In formulas,

$$\mathcal{F}u(k) := \int e^{2\pi i k\cdot x} u(x) dx \quad \text{and } \mathcal{F}^{-1}v(x) = \sum_{k\in\mathbb{Z}^d} e^{-2\pi i k\cdot x} v(k).$$

Moreover, the paraproducts $\prec, \succ$ and resonant product $\circ$ are introduced in Appendix A as well as the short-hand notations such as $\preccurlyeq = \prec + \circ$ and $\prec\succ = \prec + \succ$ which will appear on occasion.

A notation we frequently use is $f_{(\varepsilon)}$ for a function or operator to mean that the regularised object $f_\varepsilon$ and the limiting object $f$ are defined mutatis mutandis.

**Acknowledgments.**

The author acknowledges funding by the Deutsche Forschungsgemeinschaft (DFG, German Research Foundation) – CRC/TRR 388 "Rough Analysis, Stochastic Dynamics and Related Fields" – Project ID 516748464.

The author also warmly thanks Antoine Mouzard for comments on a preliminary draft of the paper.

# 2 Motivation behind the transforms

The transform $\Theta$ is a composition of three different types of transforms that are qualitatively different, sometimes applied in several different ways. We will briefly give a motivation for how they are used and where they have been successfully applied in the literature.

## 2.1 Exponential transform

This is the conceptually most simple transform and contains the renormalisation (if needed). It also neatly describes where exactly Assumption 1 comes from.

We use the formal expression

$$A(u) = (1-\Delta)u + \nabla V\cdot\nabla u + \xi u + \operatorname{div}(\rho u).$$

The main aim of this transform can be formulated as

$$\begin{array}{rl} \text{replacing } A & \text{by } e^{V+W} A e^{W} \\ \text{leads to} & \\ \text{“very singular” perturbation of } \Delta & \rightarrow \text{“less singular” perturbation of } \operatorname{div} e^{V+2W}\nabla \end{array}$$

and it is inspired by [14], [8] and [11] amongst others.

More precisely, the "inner" multiplicative term $e^W$ is a very effective and light-weight way of removing a singular linear term of the form $\Delta W\cdot u$ and, by approximately solving the elliptic KPZ type equation in Assumption 1, removing also the worst part of the quadratic term $|\nabla W|^2 u$ which appears when applying the Laplacian to $e^W$.

The second "outer" step, consisting in multiplying the differential operator $Ae^W$by $e^{V+W}$, allows us to remove the singular gradient terms $\nabla V\cdot\nabla$ and $\nabla W\cdot\nabla$, where the latter is an artifact of the first exponential transform.

In order to demonstrate this, we begin with the following simple computation which is rigorous in the smooth case. We make an ansatz $u=e^W v$ to obtain

$$\begin{aligned} A(e^W v) &= (1-\Delta)e^W v+\nabla V\cdot\nabla(e^W v)+\xi e^W v+\operatorname{div}(\rho e^W v)\\ &= e^W[(1-\Delta)v+(\nabla V+2\nabla W)\cdot\nabla v+(\xi-\Delta W-|\nabla W|^2+\nabla V\cdot\nabla W)v]+\operatorname{div}(\rho e^W v)\\ &= e^W[(1-\Delta)v+(\nabla V+2\nabla W)\cdot\nabla v+Zv]+\operatorname{div}(\rho e^W v) \qquad (13)\end{aligned}$$

where $Z$ is the remainder from the elliptic KPZ-type equation in Assumption 1.

We postpone the treatment of the term $\operatorname{div}(e^W\rho v)$ to Section 2.4 so until then we set $\rho\equiv 0$.

The linear term $e^W Z\cdot v$ in (13) formally has the regularity $-1+\delta'$ (cf. (6)) so we call it "$H^1-$ subcritical", simply because for any $\kappa>0$ there exists a constant $C(\kappa,\delta',W,Z)>0$ s.t.

$$\int v^2 e^W Z =: (v^2, e^W Z)_{B^{1-\delta'}_{1,1}} \leqslant \kappa\|v\|^2_{H^1}+C(\kappa,\delta',W,Z)\|v\|^2_{L^2}.$$

However, the gradient term is formally very irregular. In order to remedy this, we make the following observation

$$\begin{aligned} \tilde{A}(v) &= e^{W+V}A(e^W v)\\ &= e^{V+2W}[(1-\Delta)v+(\nabla V+2\nabla W)\cdot\nabla v+Zv]\\ &= (e^{V+2W}-\operatorname{div} e^{V+2W}\nabla)v+e^{V+2W}Zv\\ &= (1-\operatorname{div} e^{V+2W}\nabla)v+\tilde{Z}v,\\ &=: Lv+\tilde{Z}v \qquad (14)\end{aligned}$$

where $\tilde{Z}$ is the final object from Assumption 1.

Hence, by multiplying with the (invertible, uniformly positive) functions $e^W$, $e^{W+V}$we have reduced the study of $A$ which is a very singular (formal) perturbation of the Laplacian, to the study of $\tilde{A}$ which is a "$H^1-$subcritical" (meaning $\mathcal{C}^{-1+\delta'}$) perturbation of an elliptic operator $L$ with low-regularity coefficients (recall $e^{V+2W}\in\mathcal{C}^\delta$).

This is already a great simplification and in fact if one is interested in constructing a form domain in the case $V=0=\rho$, this is actually sufficient, see [11] which was generalised greatly in [19].

It is, however, not refined enough to precisely parametrise a domain, unless one were able to solve the elliptic KPZ-type equations in Assumption 1 with remainders of positive regularity.

## 2.2 Paracontrolled ansatz

Inspired by the theory of Paracontrolled Distributions, introduced in [10], we give a basic formulation which is sufficient in this setting. We consider the operator, formally given by

$$B(v)=\Lambda v+\tilde{Z}v$$

for $\tilde{Z}\in\mathcal{C}^{-1+\delta'}$ and $\Lambda\in\mathcal{L}(H^{1+\delta'};H^{-1+\delta'})\cap\mathcal{L}(\mathcal{D}(\Lambda);L^2)$,which has dense domain $\mathcal{D}(\Lambda)\subset H^{1+\delta}$ with inverse $\Lambda^{-1}\in\mathcal{L}(H^{-1+\delta'};H^{1+\delta'})\cap\mathcal{L}(L^2;\mathcal{D}(\Lambda))$.

We make the usual paraproduct decomposition, see the appendix,

$$\tilde{Z} v = \tilde{Z} \succ v + \tilde{Z} \preccurlyeq v$$

and if we *assume* that $v$ satisfies the ansatz

$$v - \Lambda^{-1}(\tilde{Z} \succ v) = v^\sharp \in \mathcal{D}(\Lambda) \;\text{ and }\; v \in H^{1+\delta'} \tag{15}$$

Then

$$\begin{aligned} B(v) &= \Lambda v + \tilde{Z} v \\ &= \Lambda v^\sharp + \tilde{Z} \preccurlyeq v \end{aligned}$$

and assuming (15), one has $\|B(v)\|_{L^2} \lesssim \|\Lambda v^\sharp\|_{L^2} + \|v\|_{H^1} \|\tilde{Z}\|_{\mathcal{C}^{-1+\delta'}}$, by the paraproduct estimates from Lemma 20.

In principle, as long as we make the ansatz (15) rigorous and invertible, i.e. if we have an invertible map $v = \Gamma v^\sharp$ and we take $\Lambda = L$ in (14), (having set $\rho \equiv 0$) we have formally

$$\begin{aligned} \tilde{A} \Gamma v^\sharp &= L \Gamma v^\sharp + \tilde{Z} \Gamma v^\sharp \\ &= L v^\sharp + \tilde{Z} \preccurlyeq \Gamma v^\sharp \end{aligned}$$

which would already imply

$$\|A e^W \Gamma v^\sharp\|_{L^2} \leqslant \|e^{-W-V}\|_{L^\infty} \|\tilde{A} \Gamma v^\sharp\|_{L^2} \lesssim \|L v^\sharp\|_{L^2} + \|\Gamma v^\sharp\|_{H^1} \|\tilde{Z}\|_{\mathcal{C}^{-1+\delta}}$$

and observing that for any $\lambda > 0$

$$\|\Gamma v^\sharp\|_{H^1} \leqslant c \|v^\sharp\|_{H^1} \leqslant \lambda \|L v^\sharp\|_{L^2} + \frac{c^2}{\lambda} \|v^\sharp\|_{L^2}$$

by Young's inequality and interpolation. This, together with some elliptic regularity result on $L$ would even imply

$$\|A e^W \Gamma v^\sharp\|_{L^2} + \|v^\sharp\|_{L^2} \approx \|L v^\sharp\|_{L^2} + \|v^\sharp\|_{L^2}$$

and this would already be enough to give a description of a domain of $A$ in the case $\rho \equiv 0$.

The only things that are missing is that we said that the change of variables should transform $A$ into a perturbation of the Laplacian and that we have not yet dealt with the $\operatorname{div}(\rho u)$ term.

## 2.3 Parametrix

In this section we apply a trick which is similar to Lemma 2.8 in [7] which allows us to treat a differential operator with low-regularity coefficients which are a small perturbation in $L^\infty$ of the identity as a well-behaved perturbation of the Laplacian.

To be more precise, we consider the operator $L = 1 - \operatorname{div} e^Y \nabla$ and we make the following decomposition (see the appendix for the definition and properties of the paraproducts)

$$\begin{aligned} L w &= (1 - \Delta) w - \operatorname{div}(e^Y - 1) \nabla w \\ &= (1 - \Delta) w - \operatorname{div}((e^Y - 1) \prec \nabla w) + \operatorname{div}((e^Y - 1) \succcurlyeq \nabla w) \\ &=: \Lambda w + \operatorname{div}((e^Y - 1) \succcurlyeq \nabla w) \end{aligned} \tag{16}$$

and we will see in Lemma 11 that if $\|e^Y - 1\|_{L^\infty}$ (or $\|Y\|_{L^\infty}$) is sufficiently small then we can invert $\Lambda$ on $H^\sigma$ for any $\sigma$, since the low regularity of $e^Y$ only appears in the low frequencies meaning it does not affect the regularity as it would if one were to directly invert $L$.

Observe that in (16), the term $\operatorname{div}((e^Y - 1) \succ \nabla w)$ is badly behaved, since it will always have negative regularity regardless of how good $w$ is, but we also have that

$$\|\operatorname{div}((e^Y - 1) \succcurlyeq \nabla w)\|_{H^{-1+\kappa}} \lesssim \|w\|_{H^1} \|e^Y - 1\|_{\mathcal{C}^\kappa}$$

so as long as $w$ is at least $H^1$, this term is well-defined. Importantly, the divergence is *outside* the resonant product, hence the most problematic term is actually fine.

Following this observation, we can in fact recycle the main idea from the previous section and make the paracontrolled-type ansatz

$$w = \Lambda^{-1} \operatorname{div}((e^Y - 1) \succcurlyeq \nabla w) + w^\flat \quad \text{for } w^\flat \in H^2, w \in H^{1+\kappa},$$

then (we will again be able to invert this, meaning that we have $w = \Omega w^\flat$) we have

$$L \Omega w^\flat = \Lambda w^\flat.$$

So the map $\Omega$ "replaces" the elliptic operator $L$ –which has low-regularity coefficients and hence a resolvent with limited regularisation– by the paradifferential operator $\Lambda$, where the irregular coefficients only appear in the low-frequencies.

Finally, we want to "replace" the paradifferential operator $\Lambda$ by $(1 - \Delta)$ which will be the final major step we need to achieve all the claims from Theorem 2.

This is in fact proved almost in an identical way to the invertibility of $\Lambda$. Again, observe that if $\|Y\|_{L^\infty}$ is small enough, we have that

$$\Upsilon w = w + (1 - \Delta)^{-1} \operatorname{div}((e^Y - 1) \prec \nabla w) \tag{17}$$

is invertible on $H^\sigma$ for general $\sigma$ and

$$(1 - \Delta) \Upsilon w = \Lambda w \Leftrightarrow (1 - \Delta) w = \Lambda \Upsilon^{-1} w.$$

By combining the two above transforms, we can write

$$L w = (1 - \operatorname{div} e^Y \nabla) w = (1 - \Delta) \Upsilon^{-1} \Omega^{-1} w.$$

As we shall see, these kinds of transforms allow us to transform $\tilde{A}$ and then also $A$ into a well-behaved perturbation of the Laplacian.

Note that one can completely analogously define

$$\overline{\Upsilon} w = w + (1 - \Delta)^{-1} (e^Y - 1) \prec \Delta w \tag{18}$$

which allows us to "replace" operators like $\overline{L} = 1 - e^Y \Delta$ by a Laplacian composed with some change of variables.

Observe, however, that the term $(e^Y - 1) \circ \Delta w$ is only defined when $w \in H^{2-\delta'}$, while in the previous case, $H^{1-\delta'}$ was sufficient.

## 2.4 Dealing with the div $(\rho u)$ term

Finally, we give the heuristic of how to deal with the term $\operatorname{div}(\rho v)$, where $\operatorname{div}\rho = 0$. Firstly, we give the simple observation that

$$\operatorname{div}(\rho u) = \rho \cdot \nabla u = \rho \prec \nabla u + \rho \succ \nabla u + \operatorname{div}(u \circ \rho)$$

which implies that

$$\|\operatorname{div}(\rho u)\|_{H^{-1+\delta'}} \lesssim \|u\|_{H^1} \|\rho\|_{\mathcal{C}^{-1+\delta}}$$

so it is more well-behaved than one might naively expect from a product.

However, the exponential transform interferes with this nice structure. This term behaves as follows under the exponential transform

$$\operatorname{div}(\rho e^W v) = e^W(\rho \cdot \nabla W v + \operatorname{div}(\rho v)).$$

Now while the first term can morally be treated with modifying the exponential transform, the second term actually contains a singular product between $e^W \in \mathcal{C}^\delta$ $\operatorname{div}(\rho v) \in H^{-1+\delta'}$ and that seems not to be treated in full generality without some additional tools e.g. from Paracontrolled Distributions [10].

Recalling the definition of $\tilde{A}$, we have that the div $-$ term becomes

$$\begin{aligned} e^{V+W}\operatorname{div}(\rho e^W v) &= e^{V+W}\rho \cdot \nabla(e^W v) \\ &= e^{V+2W}\nabla W \cdot \rho v + e^{V+2W}\rho \cdot \nabla v \\ &= e^{V+2W}\nabla W \cdot \rho v + e^{V+2W}\rho \prec\succ \nabla v + \operatorname{div}(e^{V+2W}\rho \circ v) + (\nabla(e^{V+2W}) \cdot \rho) \circ v \end{aligned} \tag{19}$$

leading to the bound

$$\|(19)\|_{H^{-1+\delta}} \lesssim \|v\|_{H^1}(\|\rho e^{V+2W}\|_{\mathcal{C}^{-1+\delta}} + \|e^{V+2W}\nabla W \cdot \rho\|_{\mathcal{C}^{-1+\delta'}} + \|e^{V+2W}\nabla V \cdot \rho\|_{\mathcal{C}^{-1+\delta'}}) \tag{20}$$

which is where the somewhat unnatural assumption (8) comes in which is notably satisfied when $\xi = V = 0$ and $\rho \in \mathcal{C}^{-1+\delta}$ is generic.

This bound shows us that one can simply include this term in a paracontrolled ansatz, i.e. combining it with the previous analysis, we have

$$\begin{aligned} \tilde{A}v &= Lv + \tilde{Z}v + e^{V+W}\operatorname{div}(\rho e^W v) \\ &= \Lambda v + \operatorname{div}((e^{2W+V} - 1) \succcurlyeq \nabla v) + \tilde{Z} \succ v + e^{V+W}\operatorname{div}(\rho e^W v) + \tilde{Z} \preccurlyeq v \end{aligned}$$

making the ansatz

$$v = v^\sharp + \Lambda^{-1}(\operatorname{div}((e^{2W+V} - 1) \succcurlyeq \nabla v) + \tilde{Z} \succ v + e^{V+W}\operatorname{div}(\rho e^W v))$$

we see that all the objects in the bracket are continuously bounded (depending on noise terms) from $H^1$ to $H^{-1+\delta'}$and we will see that we can make it invertible $v = \overline{\Gamma} v^\sharp$ so formally

$$\begin{aligned} \tilde{A}\overline{\Gamma}v^\sharp &= \Lambda v^\sharp + \tilde{Z} \preccurlyeq \overline{\Gamma} v^\sharp \\ &= (1 - \Delta)\Upsilon v^\sharp + \tilde{Z} \preccurlyeq \overline{\Gamma} v^\sharp \\ &\Leftrightarrow \\ \tilde{A}\overline{\Gamma}\Upsilon^{-1}v^\flat &= (1 - \Delta)v^\flat + \tilde{Z} \preccurlyeq \overline{\Gamma}\Upsilon v^\flat \end{aligned} \tag{21}$$

recalling the map $\Upsilon$ from (17).

This basically gives all the ingredients to prove the main result for $\tilde{A}$ which allows us to obtain the desired results for the actual operator $A$ by using the map $\overline{\Upsilon}$ from (18) with $Y = -W - V$. More precisely we write

$$
\begin{aligned}
A\overline{\Gamma}\Upsilon^{-1}v^{\flat} &= e^{-W-V}\tilde{A}\overline{\Gamma}\Upsilon^{-1}v^{\flat} \\
&= e^{-W-V}(1-\Delta)v^{\flat} + e^{-W-V}\tilde{Z} \preccurlyeq \overline{\Gamma}\Upsilon^{-1}v^{\flat} \\
&= (1-\Delta)\overline{\Upsilon}v^{\flat} + e^{-W-V} \succcurlyeq (1-\Delta)v^{\flat} + e^{-W-V}\tilde{Z} \preccurlyeq \overline{\Gamma}\Upsilon^{-1}v^{\flat} \\
&\Leftrightarrow \\
A\overline{\Gamma}\Upsilon^{-1}\overline{\Upsilon}^{-1}v^{\natural} &= (1-\Delta)v^{\natural} + e^{-W-V} \succcurlyeq (1-\Delta)\overline{\Upsilon}^{-1}v^{\natural} + e^{-W-V}\tilde{Z} \preccurlyeq \overline{\Gamma}\Upsilon^{-1}\overline{\Upsilon}^{-1}v^{\natural} \qquad (22)
\end{aligned}
$$

# 3 Rigorous construction of transforms

We prove all the relevant lemmas in this section which introduce rigorously the maps we introduced formally in Section 2. For $L \in \mathbb{N}$ we consider Fourier-projectors of the form

$$P_{>L} := \mathcal{F}^{-1}\mathbb{I}_{>2^L}\mathcal{F} \text{ and } P_{\leqslant L} := \mathcal{F}^{-1}\mathbb{I}_{\leqslant 2^L}\mathcal{F}.$$

**Lemma 10. (Exponential transform)** *Let $M > 0$ and $\delta > \delta' > 0$, $V_\varepsilon, V, W_\varepsilon, W$ as in Assumption 1 then*

$$e^{\pm P_{>M}V_\varepsilon} \rightarrow e^{\pm P_{>M}V}, \quad \textit{and} \quad e^{\pm P_{>M}W_\varepsilon} \rightarrow e^{\pm P_{>M}W} \quad \textit{in } \mathcal{C}^\delta \textit{ as } \varepsilon \rightarrow 0$$

*and*

$$\|e^{\pm P_{>M}Y_{(\varepsilon)}} - 1\|_{L^\infty} \lesssim 2^{-\delta' M}\|Y\|_{\mathcal{C}^\delta} \textit{ for } M \textit{ large enough depending on } Y \in \{V, W\}$$

*where importantly the $M$ can be chosen independently of $\varepsilon$.*

**Proof.** All this is quite straightforward, we omit the details. ☐

Next we show the construction and invertibility of the parametrix.

**Lemma 11.** *For all the objects and parameters in Assumption 1, there exists an $M \in \mathbb{N}$ large enough depending on the data but uniform in $\varepsilon$ s.t.*

$$\Lambda_{(\varepsilon)}w := (1-\Delta)w - \operatorname{div}\left( (e^{P_{>M}(V_{(\varepsilon)} + 2W_{(\varepsilon)})} - 1) \prec \nabla \right) w \qquad (23)$$

*and*

$$\bar{\Lambda}_{(\varepsilon)}w := (1-\Delta)w + (e^{P_{>M}(-W_{(\varepsilon)} - V_{(\varepsilon)})} - 1) \prec (1-\Delta)w \qquad (24)$$

*are invertible in $\mathcal{L}(H^\sigma; H^{\sigma-2})$ for $\sigma \in [-2, 2]$ and can be written as*

$$\Lambda_\varepsilon = (1-\Delta)\Upsilon_\varepsilon \quad \textit{and} \qquad \overline{\Lambda}_\varepsilon = (1-\Delta)\overline{\Upsilon}_\varepsilon \qquad (25)$$

*with bounded invertible maps $\Upsilon_\varepsilon, \Upsilon, \overline{\Upsilon}_\varepsilon, \overline{\Upsilon} \in \mathcal{L}(H^\sigma)$ that satisfy*

$$\Upsilon_\varepsilon \rightarrow \Upsilon \quad \textit{and} \quad \bar{\Upsilon}_\varepsilon \rightarrow \overline{\Upsilon} \quad \textit{in } \mathcal{L}(H^\sigma)$$

*and thus*

$$\Lambda_\varepsilon \rightarrow \Lambda \quad \textit{and} \quad \bar{\Lambda}_\varepsilon \rightarrow \overline{\Lambda} \quad \textit{in } \mathcal{L}(H^\sigma; H^{\sigma-2})$$

*and the analogous result holds for their bounded inverses. The same is true if one replaces the Sobolev spaces $H^\sigma$ and $H^{\sigma-2}$ by Besov spaces $B^\sigma_{p,q}$ and $B^{\sigma-2}_{p,q}$ respectively for any $p,q\in[1,\infty]$.*

**Proof.** We consider only $\Lambda$ since $\overline{\Lambda}$ is completely analogous and it is clear that one can choose $M$ large enough for both statements to hold simultaneously.

We set

$$\Upsilon_{(\varepsilon)}v := v-(1-\Delta)^{-1}\mathrm{div}\,((e^{P_{>M}(V_{(\varepsilon)}+2W_{(\varepsilon)})}-1)\prec\nabla v) \tag{26}$$

and we observe that by paraproduct estimates, Lemma 20 (see [3] to see that one can choose the implicit constant uniform in all $\sigma\in[-2,2]$)

$$\begin{aligned}\|(\Upsilon_{(\varepsilon)}-1)v\|_{H^\sigma} &\leqslant \|(e^{P_{>M}(V_{(\varepsilon)}+2W_{(\varepsilon)})}-1)\prec\nabla v\|_{H^{\sigma-1}}\\ &\leqslant C\|(e^{P_{>M}(V_{(\varepsilon)}+2W_{(\varepsilon)})}-1)\|_{L^\infty}\|v\|_{H^\sigma}\\ &\leqslant \frac{1}{2}\|v\|_{H^\sigma}\end{aligned}$$

where in the final step we have used Lemma 10 and the fact that we choose $M$ large on the data (but uniformly in $\varepsilon$). This implies that we can invert $\Upsilon_{(\varepsilon)}$ and we simply see that

$$(\Upsilon_\varepsilon-\Upsilon)v=(1-\Delta)^{-1}\mathrm{div}\,((e^{P_{>M}(V_\varepsilon+2W_\varepsilon)}-e^{P_{>M}(V+2W)})\prec\nabla v)$$

which converges strongly in $H^\sigma$. Since by construction $\Lambda_{(\varepsilon)}=(1-\Delta)\Upsilon_{(\varepsilon)}$, all the other statements follow from this. $\square$

The final technical lemma we prove is how to make the paracontrolled ansatz rigorous. The idea is motivated in Section 2.2 and is inspired by [10] and [11].

**Lemma 12. (Paracontrolled ansatz)** *For $\Lambda_{(\varepsilon)}$ as in Lemma 11 and as usual under Assumption 1, there exists an $N\in\mathbb{N}$ s.t. the maps*

$$\Phi_{(\varepsilon)}w := w-\Lambda^{-1}_{(\varepsilon)}P_{>N}\big(w\prec\tilde{Z}_{(\varepsilon)}+\mathrm{div}\,(\nabla w\preccurlyeq(e^{P_{>M}(V_{(\varepsilon)}+2W_{(\varepsilon)})}-1))+e^{P_{>M}(V_{(\varepsilon)}+W_{(\varepsilon)})}\mathrm{div}\,(\rho_{(\varepsilon)}e^{P_{>M}W_{(\varepsilon)}}v)\big)$$

*are invertible in $\mathcal{L}(H^1)$ and setting $\Gamma_{(\varepsilon)}:=\Phi^{-1}_{(\varepsilon)}$ we have that*

$$\Phi_\varepsilon\to\Phi\quad \textit{and}\ \Gamma_\varepsilon\to\Gamma\quad \textit{in}\ \mathcal{L}(H^1).$$

**Proof.** This is proved again similarly. Using the regularising properties of $\Lambda^{-1}_{(\varepsilon)}$ from Lemma 11 and the paraproduct estimates from Lemma 20 and Berstein's inequality Lemma 22, as well as the identity (19) for the div term we have

$$\begin{aligned}\|\Phi_{(\varepsilon)}(w)-w\|_{H^1} &\leqslant C\|P_{>N}(w\prec\tilde{Z}_{(\varepsilon)}+\mathrm{div}\,(\nabla w\preccurlyeq(e^{P_{>M}(V_{(\varepsilon)}+2W_{(\varepsilon)})}-1)))+\\ &\quad e^{V_{(\varepsilon)}+W_{(\varepsilon)}}\mathrm{div}\,(\rho_{(\varepsilon)}e^{W_{(\varepsilon)}}v)\|_{H^{-1}}\\ &\leqslant C2^{-\delta' N}(\|w\prec\tilde{Z}_{(\varepsilon)}\|_{H^{-1+\delta'}}+\|\nabla w\preccurlyeq(e^{P_{>M}(V_{(\varepsilon)}+2W_{(\varepsilon)})}-1)\|_{H^{\delta'}}+\\ &\quad \|e^{V_{(\varepsilon)}+W_{(\varepsilon)}}\mathrm{div}\,(\rho_{(\varepsilon)}e^{W_{(\varepsilon)}}v)\|_{H^{-1+\delta'}})\\ &\leqslant C'2^{-\delta' N}\|w\|_{H^1}(\|\tilde{Z}\|_{\mathcal{C}^{-1+\delta'}}+\|W\|_{\mathcal{C}^\delta}+\|V\|_{\mathcal{C}^\delta}+\|\rho e^{V+2W}\|_{\mathcal{C}^{-1+\delta}}+\\ &\quad \|e^{V+2W}\nabla W\cdot\rho\|_{\mathcal{C}^{-1+\delta'}}+\|e^{V+2W}\nabla V\cdot\rho\|_{\mathcal{C}^{-1+\delta'}})\\ &\leqslant \frac{1}{2}\|w\|_{H^1},\end{aligned}$$

choosing $N\in\mathbb{N}$ large enough depending on the data (and on $M$). This proves again that the maps are invertible.

It is also again straightforward to prove strong convergence

$$\begin{aligned}(\Phi_\varepsilon - \Phi)(w) &= (\Lambda_{(\varepsilon)}^{-1} - \Lambda^{-1}) P_{>N}(w \prec \tilde{Z}_\varepsilon + \operatorname{div}(\nabla w \preccurlyeq (e^{P_{>M}(V_{(\varepsilon)} + 2W_{(\varepsilon)})} - 1)) + e^{V_\varepsilon + W_\varepsilon} \operatorname{div}(\rho_\varepsilon e^{W_\varepsilon} v)) \\ &+ \Lambda^{-1} P_{>N}(w \prec (\tilde{Z}_\varepsilon - \tilde{Z}) + \operatorname{div}(\nabla w \preccurlyeq (e^{P_{>M}(V_\varepsilon + 2W_\varepsilon)} - e^{P_{>M}(V + 2W)})) + \\ & (e^{V_\varepsilon + W_\varepsilon} \operatorname{div}(\rho_\varepsilon e^{W_\varepsilon} v) - e^{V + W} \operatorname{div}(\rho e^{W} v)))\end{aligned}$$

and one sees that both terms separately converge in $H^1$ using Lemma 11 for the first term and the paraproduct estimates from Lemma 20 for the second term. The convergence of the final term follows basically by (8) in Assumption 1.

The strong convergence of $\Gamma_\varepsilon \to \Gamma$ can be proved in a completely analogous way, but in fact follows from the above. ☐

**Corollary 13.** *The operator $\Gamma$ is also invertible in $\mathcal{L}(H^\sigma; \Gamma(H^\sigma))$ for $\sigma \in [1 + \delta', 2]$ where*

$$\Gamma(H^\sigma) \subset H^{1 + \delta'} \text{ is dense but in general } \Gamma(H^\sigma) \cap H^2 = \{0\}$$

*unless the noise terms are smooth.*

**Proof.** The map $\Gamma$ is in particular injective on $H^1$, so also on $H^\sigma$.

Further, by studying

$$\Gamma(w) = w + \Lambda^{-1} P_{>N}(\Gamma(w) \prec \tilde{Z} + \operatorname{div}(\nabla \Gamma(w) \preccurlyeq (e^{P_{>M}(V + 2W)} - 1)) + e^{V + W} \operatorname{div}(\rho e^{W} \Gamma(w)))$$

one sees that no matter what the regularity of $w, \Gamma(w)$ will not be more regular than $\operatorname{reg}(\tilde{Z}) + 2$, since $\Gamma(w) \prec \tilde{Z}$ can not be more regular than $\tilde{Z}$.

The density of $H^\sigma$ in $H^{1 + \delta'}$ follows by the strong convergence of the smooth changes of variables $\Gamma_\varepsilon \to \Gamma$, similarly to [11]. ☐

# 4 Completing the proof of Theorem 2

We finally combine the heuristics from Section 2 with the rigorous formulations proven in Section 3

**Proof. (Of Theorem 2)** Note that (14) is rigorous in the regularised case, and by including the Fourier cut-off inside the exponential, we get

$$\begin{aligned}\tilde{A}_\varepsilon(v) &:= e^{P_{>M}(W_\varepsilon + V_\varepsilon)} A_\varepsilon(e^{P_{>M} W_\varepsilon} v) \\ &= (1 - \operatorname{div} e^{P_{>M}(V_\varepsilon + 2W_\varepsilon)} \nabla) v + e^{P_{>M}(W_\varepsilon + V_\varepsilon)} \nabla P_{\leqslant M} V_\varepsilon \cdot \nabla v + e^{P_{>M}(V_\varepsilon + 2W_\varepsilon)} (Z_\varepsilon + \Delta P_{\leqslant M} W_\varepsilon) v \\ &\quad + e^{P_{>M}(V_\varepsilon + W_\varepsilon)} \operatorname{div}(\rho_\varepsilon e^{P_{>M} W_\varepsilon} v) + (e^{P_{>M}(V_\varepsilon + 2W_\varepsilon)} - 1) v \\ &=: L_\varepsilon v + \tilde{V}_\varepsilon^{\leqslant M} \cdot \nabla v + \tilde{Z}_\varepsilon^M v + e^{P_{>M}(V_\varepsilon + W_\varepsilon)} \operatorname{div}(\rho_\varepsilon e^{P_{>M} W_\varepsilon} v) \end{aligned} \tag{27}$$

having redefined $\tilde{Z}_\varepsilon^M = e^{P_{>M}(V_\varepsilon + 2W_\varepsilon)} (Z_\varepsilon + \Delta P_{\leqslant M} W_\varepsilon) + e^{P_{>M}(V_\varepsilon + 2W_\varepsilon)} - 1$ and setting $\tilde{V}_\varepsilon^{\leqslant M} := e^{P_{>M}(V_\varepsilon + 2W_\varepsilon)} P_{\leqslant M} \nabla V_\varepsilon$.

Next, we use our parametrix $\Lambda_\varepsilon$ to rewrite

$$L_\varepsilon v = \Lambda_\varepsilon v + \operatorname{div}((e^{P_{>M}(V_\varepsilon + 2W_\varepsilon)} - 1) \succcurlyeq \nabla v)$$

and we assume that $v$ is paracontrolled as in Lemma 12, i.e.

$$v =: \Gamma_\varepsilon v^\sharp = v^\sharp - \Lambda_\varepsilon^{-1} P_{>N}(\Gamma_\varepsilon v^\sharp \prec \tilde{Z}_\varepsilon^M + \operatorname{div}(\nabla \Gamma_\varepsilon v^\sharp \preccurlyeq (e^{P_{>M}(V_\varepsilon + 2W_\varepsilon)} - 1)) + e^{P_{>M}(V_\varepsilon + W_\varepsilon)} \operatorname{div}(\rho_\varepsilon e^{P_{>M} W_\varepsilon} v)),$$

so (27) becomes

$$
\begin{aligned}
\tilde{A}_{\varepsilon}\Gamma_{\varepsilon}v^{\sharp} &= L_{\varepsilon}\Gamma_{\varepsilon}v^{\sharp} + \tilde{V}_{\varepsilon}^{\leqslant M}\cdot\nabla\Gamma_{\varepsilon}v^{\sharp} + \tilde{Z}_{\varepsilon}\Gamma_{\varepsilon}v^{\sharp} + e^{P_{>M}(V_{\varepsilon}+W_{\varepsilon})}\mathrm{div}\,(\rho_{\varepsilon}e^{P_{>M}W_{\varepsilon}}\Gamma_{\varepsilon}v^{\sharp})\\
&= \Lambda_{\varepsilon}\Gamma_{\varepsilon}v^{\sharp} + \mathrm{div}\,((e^{P_{>M}(V_{\varepsilon}+2W_{\varepsilon})}-1)\succcurlyeq\nabla\Gamma_{\varepsilon}v^{\sharp}) + \tilde{V}_{\varepsilon}^{M}\cdot\nabla\Gamma_{\varepsilon}v^{\sharp} + \tilde{Z}_{\varepsilon}\Gamma_{\varepsilon}v^{\sharp} +\\
&\quad e^{P_{>M}(V_{\varepsilon}+W_{\varepsilon})}\mathrm{div}\,(\rho_{\varepsilon}e^{P_{>M}W_{\varepsilon}}\Gamma_{\varepsilon}v^{\sharp})\\
&= \Lambda_{\varepsilon}v^{\sharp} + \tilde{V}_{\varepsilon}^{M}\cdot\nabla\Gamma_{\varepsilon}v^{\sharp} + \tilde{Z}_{\varepsilon}\preccurlyeq\Gamma_{\varepsilon}v^{\sharp} + P_{\leqslant N}(\mathrm{div}\,((e^{P_{>M}(V_{\varepsilon}+2W_{\varepsilon})}-1)\succcurlyeq\nabla\Gamma_{\varepsilon}v^{\sharp}) + \Gamma_{\varepsilon}v^{\sharp}\prec\tilde{Z}_{\varepsilon}^{M} +\\
&\quad \mathrm{div}\,(\nabla\Gamma_{\varepsilon}v^{\sharp}\preccurlyeq(e^{P_{>M}(V_{\varepsilon}+2W_{\varepsilon})}-1)) + e^{P_{>M}(V_{\varepsilon}+W_{\varepsilon})}\mathrm{div}\,(\rho_{\varepsilon}e^{P_{>M}W_{\varepsilon}}v))\\
&=: (1-\Delta)\Upsilon_{\varepsilon}v^{\sharp} + \overline{R}_{\varepsilon}v^{\sharp}\\
&=: (1-\Delta)v^{\flat} + R_{\varepsilon}v^{\flat}\\
&= \tilde{A}_{\varepsilon}\Gamma_{\varepsilon}\Upsilon_{\varepsilon}^{-1}v^{\flat}
\end{aligned}
$$

having set $v^{\flat} := \Upsilon_{\varepsilon}v^{\sharp}$ and recalling the identity (25) that we used in the penultimate step.

So this shows that

$$\tilde{A}_{\varepsilon}\Gamma_{\varepsilon}\Upsilon_{\varepsilon}^{-1} = (1-\Delta) + R_{\varepsilon}, \tag{28}$$

where

$$
\begin{aligned}
R_{(\varepsilon)}v^{\flat} := \; & \tilde{V}_{(\varepsilon)}^{M}\cdot\nabla\Gamma_{(\varepsilon)}\Upsilon_{(\varepsilon)}^{-1}v^{\flat} + \tilde{Z}_{(\varepsilon)}\preccurlyeq\Gamma_{(\varepsilon)}\Upsilon_{(\varepsilon)}^{-1}v^{\flat} + P_{\leqslant N}\big(\mathrm{div}\,\big((e^{P_{>M}(V_{(\varepsilon)}+2W_{(\varepsilon)})}-1)\succcurlyeq\nabla\Gamma_{(\varepsilon)}\Upsilon_{(\varepsilon)}^{-1}v^{\flat}\big) +\\
& +\Gamma_{(\varepsilon)}\Upsilon_{(\varepsilon)}^{-1}v^{\flat}\prec\tilde{Z}_{(\varepsilon)}^{M} + \mathrm{div}\,\big(\nabla\Gamma_{(\varepsilon)}\Upsilon_{(\varepsilon)}^{-1}v^{\flat}\preccurlyeq(e^{P_{>M}(V_{(\varepsilon)}+2W_{(\varepsilon)})}-1)\big) + e^{P_{>M}(V_{(\varepsilon)}+W_{(\varepsilon)})}\mathrm{div}\,(\rho_{(\varepsilon)}e^{P_{>M}W_{(\varepsilon)}}v)\big)
\end{aligned}
$$

and one obtains

$$R_{\varepsilon}\rightarrow R\in\mathcal{L}(H^{2};H^{\delta'}) \tag{29}$$

meaning that we are almost finished. Indeed, (29) can be easily proven, using the convergence of the stochastic objects, the strong convergence of the transforms proved in Section 3 and the paraproduct estimates from Lemma 20. This shows that

$$\tilde{A}_{\varepsilon}\Gamma_{\varepsilon}\Upsilon_{\varepsilon}^{-1}\rightarrow\tilde{A}\Gamma\Upsilon^{-1}\text{ in }\mathcal{L}(H^{2};L^{2})$$

and one can verify that the limit one thus obtains agrees with the rigorously defined operator

$$\tilde{A}\Gamma\Upsilon^{-1}v^{\flat} = (1-\Delta)v^{\flat} + Rv^{\flat}.$$

which allows us to define $(\tilde{A},\Gamma\Upsilon^{-1}H^{2})$ as a closed operator.

Finally, we want to transfer this result to the operator $A$ itself. Recall that

$$A_{(\varepsilon)}(u) = e^{-P_{>M}(W_{(\varepsilon)}+V_{(\varepsilon)})}\tilde{A}_{(\varepsilon)}(e^{-P_{>M}W_{(\varepsilon)}}u)$$

and

$$e^{-P_{>M}(W_{(\varepsilon)}+V_{(\varepsilon)})}(1-\Delta)w = (1-\Delta)\overline{\Upsilon}_{(\varepsilon)}w + (e^{P_{>M}(-W_{(\varepsilon)}-V_{(\varepsilon)})}-1)\succcurlyeq(1-\Delta)w$$

recalling the maps $\overline{\Upsilon}_{(\varepsilon)}$ from Lemma 11.

So, using the definition of $\tilde{A}_{\varepsilon}$ and (28) we get

$$
\begin{aligned}
A_{\varepsilon}(e^{P_{>M}W_{\varepsilon}}\Gamma_{\varepsilon}\Upsilon_{\varepsilon}^{-1}v^{\flat}) &= e^{-P_{>M}(W_{\varepsilon}+V_{\varepsilon})}\tilde{A}_{\varepsilon}\Gamma_{\varepsilon}\Upsilon_{\varepsilon}^{-1}v^{\flat}\\
&= e^{-P_{>M}(W_{\varepsilon}+V_{\varepsilon})}(1-\Delta)v^{\flat} + e^{-P_{>M}(W_{\varepsilon}+V_{\varepsilon})}R_{\varepsilon}v^{\flat}\\
&= (1-\Delta)\overline{\Upsilon}_{\varepsilon}v^{\flat} + (e^{P_{>M}(-W_{\varepsilon}-V_{\varepsilon})}-1)\succcurlyeq(1-\Delta)v^{\flat} + R'_{\varepsilon}v^{\flat}
\end{aligned}
$$

for $R'_{(\varepsilon)} = e^{-P_{>M}(W_{(\varepsilon)}+V_{(\varepsilon)})}R_{(\varepsilon)}$.

Finally, after setting $v^\natural = \bar{\Upsilon}_\varepsilon v^\flat$ and defining

$$\Theta_{(\varepsilon)} := e^{P_{>M}(W_{(\varepsilon)} + W'_{(\varepsilon)})} \Gamma_{(\varepsilon)} \Upsilon_{(\varepsilon)}^{-1} \bar{\Upsilon}_{(\varepsilon)}^{-1} \tag{30}$$

we have

$$A_{(\varepsilon)} \Theta_{(\varepsilon)} v^\natural = \; (1-\Delta) v^\natural + (e^{P_{>M}(-W_{(\varepsilon)} - V_{(\varepsilon)})} - 1) \succcurlyeq (1-\Delta) \bar{\Upsilon}_{(\varepsilon)}^{-1} v^\natural + R'_{(\varepsilon)} \bar{\Upsilon}_{(\varepsilon)}^{-1} v^\natural$$

and as before, one readily checks that $(A_\varepsilon \Theta_\varepsilon - (1-\Delta)) \rightarrow (A\Theta - (1-\Delta))$ in $\mathcal{L}(H^2; H^\delta)$ which, together with the invertibility of $\Theta$, implies the norm resolvent convergence $A_\varepsilon \rightarrow A$.

Moreover, the above result immediately implies

$$\|A\Theta v^\natural\|_{L^2} + \|\Theta v^\natural\|_{L^2} \approx \|v^\natural\|_{H^2}$$

and from the invertibility of $\Theta$, one can straightforwardly prove that $\Theta(H^2)$ is dense in $L^2$ as in [11].

This finishes the proof. □

# 5 Some examples

We give some fairly general examples of potentials that satisfy Assumption 1. Firstly we consider generic potentials i.e. where there is no special (e.g. stochastic) structure and no relation between the different terms. Then we consider the case of generalised "subcritical noises" as constructed in [19].

## 5.1 Generic potentials

If we make the following assumptions then the structural assumptions from Assumption 1 are automatically satisfied. For simplicity we restrict ourselves to $d \geqslant 2$ with some straightforward extensions to $d = 1$. We assume that the potential terms satisfy either of the following assumptions:

**Assumption (I):**

$$\begin{aligned} \xi \in B_{p,\infty}^{-1+\delta}(\mathbb{T}^d), V \in B_{q,\infty}^{1-\delta'}(\mathbb{T}^d), \;& \rho \in B_{r,\infty}^{-\delta'}(\mathbb{T}^d; \mathbb{R}^d) \\ \text{where} \;& 0 < \delta' < \delta < 1, \\ p, r > d, \;& \frac{1}{q} < \frac{1-\delta-\delta'}{d}, \\ \text{and} \;& \\ \frac{1}{s} := \frac{1}{p} + \frac{1}{r} < \frac{1}{d}, \qquad & \frac{1}{\tilde{s}} := \frac{1}{q} + \frac{1}{r} < \frac{1}{d}. \end{aligned}$$

**Assumption (II):**

$$\begin{aligned} \xi \in \mathcal{C}^{-\frac{1}{2}+\delta}(\mathbb{T}^d), V \in \mathcal{C}^{\frac{1}{2}+\delta'}(\mathbb{T}^d), \;& \rho \in B_{r,\infty}^{-\frac{1}{2}-\delta''}(\mathbb{T}^d; \mathbb{R}^d) \\ \text{where} \;& 0 < \delta'' < \delta' < \delta < 1 \\ \text{and} \;& \frac{d}{r} = \frac{1}{2} - \delta' - \delta'' \end{aligned}$$

**Remark 14.** In Assumption (I), we have a comparatively weak assumption on $\xi$ but strong assumptions on $V, \rho$. Assumption (II), on the other hand, assumes $V$ to be less regular but assumes $\xi$ to be better.

Clearly one could assume $\xi, V$ to be regular enough and $\rho$ to be only $\mathcal{C}^{-1+\delta'}$, but that case is simpler.

We will show that under these asumptions, we can actually find an exact solution to the elliptic KPZ-type equation

$$(\lambda - \Delta) W - |\nabla W|^2 + \nabla W \cdot \nabla V + \xi = 0 \qquad (31)$$

for sufficiently large $\lambda > 0$ in $B_{p,\infty}^{1+\delta} \hookrightarrow \mathcal{C}^{\delta}$ under Assumption (I) or $W \in \mathcal{C}^{\frac{3}{2}+\delta}$ assuming (II). In the former case, we have

$$\|e^{V+2W} \nabla W \cdot \rho\|_{B_{s,\infty}^{-\delta'}} + \|e^{V+2W} \nabla V \cdot \rho\|_{B_{\tilde{s},\infty}^{-\delta'}} \;\lesssim\; \|e^{V+2W}\|_{\mathcal{C}^{\delta}} \|\rho\|_{B_{r,\infty}^{-\delta'}} \big( \|\xi\|_{B_{p,\infty}^{-1+\delta}} + \|V\|_{B_{q,\infty}^{1+\delta}} \big)$$

which is the final point (8) in Assumption 1, the point (7) being analogous. Checking Assumption 1 under (II) is simpler.

**Lemma 15.** *There exists a large $\lambda > 0$ depending on the data s.t. (31) has a unique solution*

- $W \in B_{p,\infty}^{1+\delta}$ *assuming (I),*
- $W \in \mathcal{C}^{\frac{3}{2}+\delta}$ *assuming (II).*

**Proof.** We set up a fixed point and define the map

$$\Phi(W) := -(\lambda - \Delta)^{-1} (|\nabla W|^2 + \nabla W \cdot \nabla V + \xi)$$

**Part (I):** Under the Assumption (I), setting $\frac{1}{\nu} = \frac{1}{p} + \frac{1}{q}$, we use the simple bound

$$\|(\lambda - \Delta)^{-1} f\|_{B_{\mu,\infty}^{\beta}} \lesssim \lambda^{-\kappa} \|f\|_{B_{\mu,\infty}^{\beta-2+\kappa}} \quad \kappa < \beta, \quad \mu \in [1, \infty]. \qquad (32)$$

Then there exists a small $\alpha > 0$ s.t.

$$\begin{aligned} \|\Phi(W)\|_{B_{p,\infty}^{1+\delta}} &\leqslant C \|\xi\|_{B_{p,\infty}^{-1+\delta}} + C\lambda^{-\alpha} \||\nabla W|^2\|_{B_{\frac{p}{2},\infty}^{\delta}} + C\lambda^{-\alpha} \|\nabla W \cdot \nabla V\|_{B_{\nu,\infty}^{\delta}} \\ &\leqslant C \|\xi\|_{B_{p,\infty}^{-1+\delta}} + C\lambda^{-\alpha} \|\nabla W\|_{B_{p,\infty}^{\delta}}^2 + C\lambda^{-\alpha} \|V\|_{B_{q,\infty}^{1-\delta'}} \|\nabla W\|_{B_{p,\infty}^{\delta}} \end{aligned}$$

using Besov embedding and product inequalities, and

$$\|\Phi(W) - \Phi(W')\|_{B_{p,\infty}^{1+\delta}} \leqslant C\lambda^{-\alpha} \|\nabla(W - W')\|_{B_{p,\infty}^{\delta}} \|\nabla(W + W')\|_{B_{p,\infty}^{\delta}} + C\lambda^{-\alpha} \|V\|_{B_{q,\infty}^{1-\delta'}} \|\nabla(W - W')\|_{B_{p,\infty}^{\delta}}$$

so it's clear that for $\lambda > 0$ large enough depending on the data that we can get a fixed point.

**Part (II):** This is simpler. Under Assumption (II) we have

$$\begin{aligned} \|\Phi(W)\|_{\mathcal{C}^{\frac{3}{2}+\delta}} &\leqslant C \|\xi\|_{\mathcal{C}^{-\frac{1}{2}+\delta}} + C\lambda^{-\alpha} \||\nabla W|^2\|_{\mathcal{C}^{\frac{1}{2}+\delta}} + C\lambda^{-\alpha} \|\nabla W \cdot \nabla V\|_{\mathcal{C}^{-\frac{1}{2}+\delta}} \\ &\leqslant C \|\xi\|_{\mathcal{C}^{-\frac{1}{2}+\delta}} + C\lambda^{-\alpha} \|\nabla W\|_{\mathcal{C}^{\frac{1}{2}+\delta}}^2 + C\lambda^{-\alpha} \|\nabla V\|_{\mathcal{C}^{-\frac{1}{2}+\delta}} \|\nabla W\|_{\mathcal{C}^{\frac{1}{2}+\delta}} \end{aligned}$$

and

$$\|\Phi(W) - \Phi(W')\|_{\mathcal{C}^{\frac{3}{2}+\delta}} \leqslant C\lambda^{-\alpha} \|\nabla(W - W')\|_{\mathcal{C}^{\frac{1}{2}+\delta}} \|\nabla(W + W')\|_{\mathcal{C}^{\frac{1}{2}+\delta}} + C\lambda^{-\alpha} \|\nabla V\|_{\mathcal{C}^{-\frac{1}{2}+\delta}} \|\nabla(W - W')\|_{\mathcal{C}^{\frac{1}{2}+\delta}}$$

using the standard product bounds for Hölder spaces.

Then one readily gets a fixed point for $\lambda > 0$ sufficiently large. □

### 5.2 Generalised Anderson Hamiltonians

In this section, we set $V=0=\rho$ but the construction also works if they are regular enough.

The study of singular operators of the form $-\Delta+\xi$ for $\xi$ being spatial white noise on $\mathbb{T}^2$ including a rigorous construction of a domain was initiated in [2] and used the theory of Paracontrolled Distributions [10]. This was reformulated and extended to $\xi$ being spatial white noise on $\mathbb{T}^3$ in [11] which is more involved and includes the trick of combining the exponential transform with a paracontrolled ansatz which also appears in the current article. It was also first observed in that article that by using only an exponential transform, one can construct the operator via its form domain, a result that was extended to the full sub-critical regime of generalised noises in [19].

In fact, in that paper, the authors construct solutions of the elliptic KPZ-type equation (2) with $V=0=\rho$ and $\xi_\varepsilon$ a regularisation of a *subcritical noise*, see Section 3.1 of that paper. The authors also construct the (in general diverging) renormalisation constants $c_\varepsilon$ and show that

$$\mathcal{H}_{\xi_\varepsilon} := 1-\Delta+\xi_\varepsilon+c_\varepsilon \to \mathcal{H}_\xi$$

in the norm resolvent sense for a self-adjoint limit $\mathcal{H}_\xi$ defined as a Friedrichs extension of a quadratic form.

One result of the current paper is that (using their result on the convergence of the noise terms), we can construct a refined change of variables $\Theta_{(\varepsilon)}$ so we actually have the stronger convergence

$$\mathcal{H}_{\xi_\varepsilon}\Theta_\varepsilon \to \mathcal{H}_\xi\Theta \text{ in } \mathcal{L}(H^2;L^2)$$

which also allows to parametrise a domain for the limiting operator.

**Remark 16.** Probably it is achievable but somewhat technical to prove a general convergence result for subcritical stochastic terms $V,\rho\neq 0$ however one loses the self-adjointness of the operator in the case $\rho\neq 0$. In fact, from the usual "power counting" arguments from [13] and subsequent papers, one would expect Assumption 1 to be true for a rather generic subcritical noise $V$, but because of (8), we need there to be some cancellation in the scalar products $e^{V+2W}\nabla V\cdot\rho$ and $e^{V+2W}\nabla W\cdot\rho$, so $\rho$ should depend on $W,V$ in some sense. Alternatively we could "only" assume $\rho\in\mathcal{C}^{-\delta''}(\mathbb{T}^d;\mathbb{R}^d)$, for $\delta''<\delta'$ small enough, which by power counting would suggest (8) to hold.

# 6 Selfadjointness and applications to singular NLS

In this section, we show how to apply the methods from the previous sections to nonlinear Schrödinger equations (NLS) with singular potentials

$$i\partial_t u + Au = -|u|^2 u \;\text{ on } \mathbb{T}^d$$

for

$$A(u) \text{``="} (1-\Delta)u + \nabla V\,\nabla u + \xi u$$

with $V\in\mathcal{C}^\delta,\xi\in\mathcal{C}^{-2+\delta}$ and $\rho=0$ which under Assumption 1 we have rigorously constructed.

**Proof. (of Theorem 6)** The exponential transform leads us to

$$e^{V+W}A(e^W v) \;=(1-\operatorname{div} e^{2W+V}\nabla)v+\tilde{Z}v$$

which we can reframe as saying that the operator $e^W A e^W$ (and also $A$) is symmetric on $L^2(e^V dx)=:$ $L^2_V$ while it is *not* symmetric on $L^2(dx)$.

This can be seen by the computation

$$\begin{aligned}\langle e^W A e^W v, w \rangle_V &:= (e^V e^W A e^W v, w)_{L^2}\\ &= ((1 - \operatorname{div} e^{2W+V} \nabla) v + \tilde{Z} v, w)_{L^2}\\ &= (v, (1 - \operatorname{div} e^{2W+V} \nabla) w + \tilde{Z} w)_{L^2}\\ &= \langle v, e^W A e^W w \rangle_V\end{aligned}$$

which shows symmetry of the operator. The existence of a unique self-adjointness operator can be shown using the fact that the quadratic form

$$\langle e^W A e^W v, v \rangle_V = \int e^{V+2W} |\nabla v|^2 + (\tilde{Z} v, v)$$

is coercive. In fact one has for sufficiently large $C > 0$

$$C \|v\|_{L_V^2}^2 + \langle e^W A e^W v, v \rangle_V \geqslant (\|e^{-2W}\|_{L^\infty} \|e^{-V}\|_{L^\infty} + 1)^{-1} \|v\|_{H^1}^2 \tag{33}$$

using that the other term is bounded by $|(\tilde{Z} v, v)| \lesssim \|\tilde{Z}\|_{\mathcal{C}^{-1+\delta}} \|v^2\|_{B_{1,1}^{1-\delta}} \lesssim \|v\|_{H^{1-\delta'}}^2$ and is thus lower order w.r.t. the $H^1$ norm. Then the existence of a unique self-adjoint extension follows by applying e.g. Theorem VIII.15 in [26]. □

We can describe the space $\operatorname{dom}_V(A) = \{f : A f \in L_V^2\}$ with the equivalent norm

$$\|A f\|_{L_V^2} = \|e^V A f\|_{L^2} \approx \|e^V e^W A f\|_{L^2}.$$

Setting $f = e^W g$ we have as in Section 2.1

$$\|e^V e^W A f\|_{L^2} = \|e^V e^W A e^W g\|_{L^2} = \|L g + \tilde{Z} g\|_{L^2}$$

having set $L := (1 - \operatorname{div} e^{2W+V} \nabla) v$ and using the paracontrolled ansatz

$$g = L^{-1}(g \prec \tilde{Z}) + g^\sharp \tag{34}$$

which, after suitable truncation as in Section 3, we can write as $g = \Gamma g^\sharp$ for an invertible map $\Gamma$. This leads to

$$\|L g + \tilde{Z} g\|_{L^2} = \|L \Gamma g^\sharp + \tilde{Z} \Gamma g^\sharp\|_{L^2} = \|L g^\sharp + \tilde{Z} \preccurlyeq \Gamma g^\sharp\|_{L^2} \lesssim \|L g^\sharp\|_{L^2} + \|g^\sharp\|_{L^2}$$

so

$$\|A e^W \Gamma g^\sharp\|_{L_V^2} \lesssim \|L g^\sharp\|_{L^2} + \|g^\sharp\|_{L^2}$$

and one can prove also the reverse inequality

$$\|L g^\sharp\|_{L^2} \lesssim \|A e^W \Gamma g^\sharp\|_{L_V^2} + \|g^\sharp\|_{L^2}$$

in the same way.

Further, we know that we can give a parametrisation of $\operatorname{dom}_V L$ using the parametrix i.e. there is a map $\Upsilon$ s.t. $\operatorname{dom}_V L = \operatorname{dom} L = \Upsilon H^2$, see Section 2.3. We need an additional embedding result on the domain of $L$.

**Lemma 17. ($L^\infty$ embedding)** *For any $u \in \operatorname{dom}(L)$ on $\mathbb{T}^d$ with $d \leqslant 3$ we have*

$$\|u\|_{L^\infty} \lesssim \|L u\|_{L^2} + \|u\|_{L^2}.$$

*Moreover, the map* $\Gamma$ *implicitly defined by (34) i.e.*

$$\Gamma g^{\sharp} = L^{-1} P_{>N}(\Gamma g^{\sharp} \prec \tilde{Z}) + g^{\sharp}$$

*satisfies the bound*

$$\|\Gamma^{\pm 1} g^{\sharp}\|_{L^{\infty}} \lesssim \|g^{\sharp}\|_{L^{\infty}}. \tag{35}$$

**Proof.** As in Section 2.3, we rewrite $L = \Lambda + \mathrm{div}\,(e^{2W+Y} \succcurlyeq \nabla)$ and by Lemma 11 after truncation (which we omit to lighten the notation) $\Lambda$ can be inverted on any $B^s_{p,q}$. Thus setting

$$u = \Lambda^{-1} \mathrm{div}\,(e^{2W+Y} \succcurlyeq \nabla u) + u^{\sharp}$$

there exists an inverse $\Omega g^{\sharp} = g$ (again after a harmless truncation) which is bounded on $B^1_{p,p}$ for any $p$. In this notation

$$L = \Lambda \Omega$$

so by Besov embedding, Lemma 23,

$$\|u\|_{L^{\infty}} \lesssim \|u\|_{B^1_{d+\kappa,d+\kappa}} \lesssim \|\Omega u\|_{B^1_{d+\kappa,d+\kappa}} \lesssim \|\Lambda \Omega u\|_{B^{-1}_{d+\kappa,d+\kappa}} \lesssim \|L g\|_{L^2} + \|g\|_{L^2},$$

as claimed.

Using the same argument, we can also show a Schauder estimate

$$\|L^{-1} f\|_{\mathcal{C}^{1+\delta}} \lesssim \|f\|_{\mathcal{C}^{-1+\delta}}. \tag{36}$$

The bound on $\Gamma$ follows from this, since by Lemma 20

$$\|\Gamma g^{\sharp}\|_{L^{\infty}} \leqslant \|g^{\sharp}\|_{L^{\infty}} + 2^{-N\delta'} \|\Gamma g^{\sharp} \prec \tilde{Z}\|_{\mathcal{C}^{-1+\delta'}} \leqslant C \|g^{\sharp}\|_{L^{\infty}} + \frac{1}{2} \|\Gamma g^{\sharp}\|_{L^{\infty}}$$

for $N$ chosen depending on the $\|\tilde{Z}\|_{\mathcal{C}^{-1+\delta'}}$. □

This immediately implies the same bound for the domain of $A$.

**Corollary 18.** *In the setting above, we have*

$$\|u\|_{L^{\infty}} \lesssim \|A u\|_{L^{\infty}} + \|u\|_{L^2}.$$

**Proof.** By boundedness of $e^{-W}$ and $L^{\infty}$ boundedness of $\Gamma^{-1}$

$$\|u\|_{L^{\infty}} \lesssim \|e^W \Gamma u\|_{L^{\infty}} \lesssim \|L e^W \Gamma u\|_{L^2} + \|u\|_{L^2} \lesssim \|A u\|_{L^2} + \|u\|_{L^2}.$$ □

Now we turn our attention to non-linear Schrödinger type equations driven by $A$ i.e. equations like

$$i \partial_t u + A u = - |u|^2 u \tag{37}$$

although one can consider more general nonlinearities $f(u)$ with relatively minimal assumptions.

Since $A$ is self-adjoint on $L^2_V$, we can define the unitary group $e^{itA}$ on that space, so

$$u = e^{itA} u_0 \Leftrightarrow \begin{cases} (i \partial_t + A) u & = 0 \\ u(0) & = u_0 \end{cases}.$$

Moreover, we can write a mild formulation for (37) as

$$u(t) = e^{itA}u_0 + i\int_0^t e^{i(t-s)A}u|u|^2(s)ds \tag{38}$$

where the equality holds in $L^2_V$ and the solution will be constructed as a fixed point in a suitable space.

Moreover, we want to say that up to a short time $T > 0,$ a solution $u$ actually exists and is the limit of

$$u_\varepsilon(t) = e^{itA_\varepsilon}u_0^\varepsilon + i\int_0^t e^{i(t-s)A_\varepsilon}u_\varepsilon|u_\varepsilon|^2(s)ds \tag{39}$$

with suitably chosen initial data $u_0^\varepsilon$.

Note that the operator $A_\varepsilon$ is self-adjoint on $L^2_{V_\varepsilon}$ but not on $L^2$, hence the mild formulation holds on a space depending on $\varepsilon$. However, all the $L^2$ spaces have equivalent norm albeit with quite different scalar products. As a preliminary step, we prove that the unitary group $e^{itA_\varepsilon}$ converges to the unitary group $e^{itA}$ strongly which would be automatic from Theorem VIII.21 [26] but we need a simple additional argument since the groups are unitary with respect to different scalar products.

**Lemma 19.** *We have that*

$$e^{itA_\varepsilon}g \to e^{itA}g \text{ in } L^2$$

*for all $t \in \mathbb{R}$ and all $g \in L^2$. One has that $e^{itA_\varepsilon}$ and $e^{itA}$ are strongly continuous groups on $L^2$ (one loses unitarity since we change the scalar product).*

**Proof.** By self-adjointness and Stone's theorem, we have that $e^{itA_\varepsilon}$ is self-adjoint on $L^2_{V_\varepsilon}$ which implies that

$$e^{V_\varepsilon}e^{itA_\varepsilon}e^{-V_\varepsilon} = e^{ite^{V_\varepsilon}A_\varepsilon e^{-V_\varepsilon}}$$

is self-adjoint on $L^2$. In fact, one computes

$$(e^{V_\varepsilon}e^{itA_\varepsilon}e^{-V_\varepsilon}u, e^{V_\varepsilon}e^{-V_\varepsilon}v) = \langle e^{itA_\varepsilon}e^{-V_\varepsilon}u, e^{-V_\varepsilon}v\rangle_{V_\varepsilon} = \langle e^{-V_\varepsilon}u, e^{itA_\varepsilon}e^{-V_\varepsilon}v\rangle_{V_\varepsilon} = (u, e^{V_\varepsilon}e^{itA_\varepsilon}e^{-V_\varepsilon}v)$$

so it is symmetric and bounded and hence self-adjoint on $L^2$.

Theorem 2 implies $A_\varepsilon^{-1} \to A^{-1}$ in $\mathcal{L}(L^2)$ and we claim that this implies also that the self-adjoint operators

$$(e^{V_\varepsilon}A_\varepsilon e^{-V_\varepsilon})^{-1} \to (e^V A e^{-V})^{-1} \text{ converge in } \mathcal{L}(L^2).$$

In fact, this follows simply by writing the difference

$$\begin{aligned}(e^V A e^{-V})^{-1} - (e^{V_\varepsilon}A_\varepsilon e^{-V_\varepsilon})^{-1} &= e^V A^{-1}e^{-V} - e^{V_\varepsilon}A_\varepsilon^{-1}e^{-V_\varepsilon}\\ &= (e^V A^{-1} - e^{V_\varepsilon}A_\varepsilon^{-1})\,e^{-V} - e^{V_\varepsilon}A_\varepsilon^{-1}(e^{-V_\varepsilon} - e^{-V})\\ &= ((e^V - e^{V_\varepsilon})A^{-1} + e^{V_\varepsilon}(A^{-1} - A_\varepsilon^{-1}))\,e^{-V} - e^{V_\varepsilon}A_\varepsilon^{-1}(e^{-V_\varepsilon} - e^{-V})\end{aligned}$$

and observing that every term goes to zero in $\mathcal{L}(L^2)$.

Now, standard functional calculus, see e.g. Theorem VIII.21 [26], implies that the unitary groups

$$e^{V_\varepsilon}e^{itA_\varepsilon}e^{-V_\varepsilon} = e^{ite^{V_\varepsilon}A_\varepsilon e^{-V_\varepsilon}} \to e^{ite^V A e^{-V}} = e^V e^{itA}e^{-V} \text{ converge strongly in } L^2.$$

Finally, undoing the multiplication as above, we can conclude that also $e^{itA_\varepsilon} \to e^{itA}$ strongly in $L^2$.

The unitarity clearly is false if one changes the scalar product, but the group property is preserved since

$$e^{isA}e^{itA} = e^{-V}e^{V}e^{isA}e^{-V}e^{V}e^{itA}e^{-V}e^{V} = e^{-V}e^{V}e^{i(s+t)A}e^{-V}e^{V} = e^{i(s+t)A}. \qquad \square$$

This allows us to prove the local well-posedness of (37) as in [11].

**Proof. (of Theorem 8)** We prove that there exists a unique solution to (38) in the space

$$C_T\,\mathrm{dom}_V A \cap C_T^1 L^2$$

Writing

$$\Psi u(t) := e^{itA}u_0 + i\int_0^t e^{isA}u|u|^2(t-s)ds, \tag{40}$$

we make a manipulation of the inhomogeneous term using integration by parts in time (it is not possible to apply the operator $A$ to the nonlinearity but the time derivative is fine)

$$\begin{aligned} i\int_0^t e^{isA}u|u|^2(t-s)ds = & -A^{-1}\int_0^t \partial_s(e^{isA})u|u|^2(t-s)ds \\ = & -A^{-1}\int_0^t e^{isA}\partial_s(u|u|^2)(t-s)ds - A^{-1}(e^{itA}u_0^3 - u|u|^2(t)). \end{aligned}$$

Thus we have straightforwardly using the $L^2$ boundedness of $e^{isA}$ and that it commutes with $A$

$$\|\Psi u\|_{C_T\mathcal{D}_V(A)\cap C_T^1L_V^2} \lesssim \|Au_0\|_{L^2} + \|u_0\|_{L^6}^3 + \|u(t)\|_{L^6}^3 + \int_0^t \|u\|_{L^\infty}^2\|\partial_t u\|_{L^2}.$$

The $C_T^1L_V^2$ bound follows more directly from (40) by applying the time derivative.

To prove the contraction property, we need to control

$$u|u|^2(t) - v|v|^2(t) = \int_0^t \partial_s(u|u|^2(s) - v|v|^2(s))$$

in addition to the usual terms involving the time integral. One thus gets

$$\begin{aligned} \|\Psi u - \Psi v\|_{C_T\mathcal{D}_V(A)\cap C_T^1L_V^2} \lesssim & T(\|\partial_t u - \partial_t v\|_{L^\infty L^2}(\|u\|_{L^\infty L^\infty}^2 + \|v\|_{L^\infty L^\infty}^2) + \|u-v\|_{L^\infty L^\infty}(\|v\|_{L^\infty L^\infty} + \\ & \|u\|_{L^\infty L^\infty})(\|\partial_t u\|_{L^\infty L^2} + \|\partial_t v\|_{L^\infty L^2})) \\ \lesssim & T\|\Psi u - \Psi v\|_{C_T\mathcal{D}_V(A)\cap C_T^1L_V^2}(\|u\|_{C_T\mathcal{D}_V(A)\cap C_T^1L_V^2}^2 + \|v\|_{C_T\mathcal{D}_V(A)\cap C_T^1L_V^2}^2) \end{aligned}$$

having used the $L^\infty$ embedding from Lemma 17.

As in [11] or [7] one can choose $M>0$ large and then $T>0$ small depending on $\|Au_0\|_{L^2}$ and the other parameters to prove that $\Psi$ is indeed a contraction on the space $C_T\mathcal{D}_V(A)\cap C_T^1L_V^2$. The time continuity follows by Stone's theorem and the self-adjointness of $A$.

This gives a local in time analytically strong solution to (37). We can further prove that it is the limit of the solutions of smoothened, classically well-posed, equations.

By choosing the initial data as $u_{0,\varepsilon} := A_\varepsilon^{-1} A u_0$, we can achieve analogously to Theorem 3.6 in[11]

$$\begin{cases} u_\varepsilon & \to & u & \text{in } C_T L^2 \\ A_\varepsilon u_\varepsilon & \to & Au & \text{in } C_T L^2 \\ \partial_t u_\varepsilon & \to & \partial_t u & \text{in } C_T L^2 \end{cases}.$$

This follows by using the strong convergence the groups proved in Lemma 19 and Gronwall, we omit the details and refer to the reader to the proof of Theorem 3.6 in[11]. □

# A Besov spaces and paraproducts

We collect some elementary results about paraproducts and Besov spaces and refer the readers for example to [3] or [10] for more details on related topics and especially the application to singular SPDEs.

We first recall the definition of Littlewood-Paley blocks and denoted by $\chi$ and $\rho$ two non-negative smooth and compactly supported radial functions $\mathbb{R}^d \to \mathbb{C}$ such that

1. The support of $\chi$ is contained in a ball and the support of $\rho$ is contained in an annulus $\{x \in \mathbb{R}^d : a \leqslant |x| \leqslant b\}$
2. For all $\xi \in \mathbb{R}^d$, $\chi(\xi) + \sum_{j \geqslant 0} \rho(2^{-j}\xi) = 1$;
3. For $j \geqslant 1$, $\chi(\cdot)\rho(2^{-j}\cdot) = 0$ and $\rho(2^{-j}\cdot)\rho(2^{-i}\cdot) = 0$ for $|i-j| > 1$.

The Littlewood-Paley blocks $(\Delta_j)_{j \geqslant -1}$ associated to $f \in \mathcal{S}'(\mathbb{T}^d)$ are defined by $\Delta_{-1} f := \mathcal{F}^{-1}(\chi \mathcal{F} f)$ and $\Delta_j f := \mathcal{F}^{-1}(\rho(2^{-j}\cdot)\mathcal{F} f)$ for $j \geqslant 0$. We also set, for $f \in \mathcal{S}'(\mathbb{T}^d)$ and $j \geqslant -1$ that $S_j f := \sum_{i=-1}^{j-1} \Delta_i f$. Then the Besov space with parameters $p, q \in [1, \infty)$, $\alpha \in \mathbb{R}$ can now be defined as $B_{p,q}^\alpha(\mathbb{T}^d) := \{u \in \mathcal{S}'(\mathbb{T}^d) : \|u\|_{B_{p,q}^\alpha} < \infty\}$, where the norm is defined as

$$\|u\|_{B_{p,q}^\alpha} := \left( \sum_{k \geqslant -1} (2^{\alpha k} \|\Delta_k u\|_{L^p})^q \right)^{\frac{1}{q}},$$

with the obvious modification for $q = \infty$. We also define the *Besov-Hölder* spaces $\mathcal{C}^\alpha := B_{\infty,\infty}^\alpha$, which for $\alpha \in (0,1)$ agree with the usual Hölder spaces $C^\alpha$. Using this notation, we can formally decompose the product $f \cdot g$ of two distributions $f$ and $g$ as

$$f \cdot g = f \prec g + f \circ g + f \succ g,$$

where

$$f \prec g := \sum_{j \geqslant -1} S_{j-1} f \Delta_j g \quad \text{and} \quad f \succ g := \sum_{j \geqslant -1} \Delta_j f S_{j-1} g$$

are referred to as the *paraproducts*, whereas

$$f \circ g := \sum_{j \geqslant -1} \sum_{|i-j| \leqslant 1} \Delta_i f \Delta_j g$$

is called the *resonant product*. An important point is that the paraproduct terms are always well defined whatever the regularity of $f$ and $g$. The resonant product, on the other hand, is a priori only well defined if the sum of their regularities is positive. We also, on occasion, write

$$\preccurlyeq := \prec + \circ \text{ and } \prec\succ := \prec + \succ. \tag{41}$$

We collect some results.

**Lemma 20.** **( [20], Theorem 3.17)** *Let $\alpha, \alpha_1, \alpha_2 \in \mathbb{R}$ and $p, p_1, p_2, q \in [1, \infty]$ be such that*

$$\alpha_1 \neq 0 \quad \alpha = (\alpha_1 \wedge 0) + \alpha_2 \quad \text{and} \quad \frac{1}{p} = \frac{1}{p_1} + \frac{1}{p_2}.$$

*Then we have the bound*

$$\| f \prec g \|_{B^{\alpha}_{p,q}} \lesssim \| f \|_{B^{\alpha_1}_{p_1,\infty}} \| g \|_{B^{\alpha_2}_{p_2,q}}$$

*and in the case where $\alpha_1 + \alpha_2 > 0$ we have the bound*

$$\| f \circ g \|_{B^{\alpha_1 + \alpha_2}_{p,q}} \lesssim \| f \|_{B^{\alpha_1}_{p_1,\infty}} \| g \|_{B^{\alpha_2}_{p_2,q}}.$$

**Remark 21.** In this paper, we mainly consider the *Bessel potential spaces,* $H^s = B^s_{2,2}$ and *Hölder Besov* spaces $\mathcal{C}^s = B^s_{\infty,\infty}$; also $B^s_{p,\infty}$ in Section 5.1.

**Lemma 22.** **(Bernstein Inequality)** *Let $\mathcal{A}$ be an annulus and $\mathcal{B}$ be a ball. For any $k \in \mathbb{N}, \lambda > 0$, and $1 \leqslant p \leqslant q \leqslant \infty$ we have*

1. *if $u \in L^p(\mathbb{T}^d)$ is such that $\operatorname{supp}(\mathcal{F}u) \subset \lambda \mathcal{B}$ then*

$$\max_{\mu \in \mathbb{N}^d : |\mu| = k} \| \partial^{\mu} u \|_{L^q} \lesssim_k \lambda^{k + d\left(\frac{1}{p} - \frac{1}{q}\right)} \| u \|_{L^p}$$

2. *if $u \in L^p(\mathbb{T}^d)$ is such that $\operatorname{supp}(\mathcal{F}u) \subset \lambda \mathcal{A}$ then*

$$\lambda^k \| u \|_{L^p} \lesssim_k \max_{\mu \in \mathbb{N}^d : |\mu| = k} \| \partial^{\mu} u \|_{L^p}.$$

**Lemma 23.** **(Besov Embedding)** *Let $\alpha < \beta \in \mathbb{R}$, $q_1 \leqslant q_2$, and $p > r \in [1, \infty]$ be such that $\beta = \alpha + d(1/r - 1/p)$, then we have the following bound*

$$\| f \|_{B^{\alpha}_{p,q_2}(\mathbb{T}^d)} \lesssim \| f \|_{B^{\beta}_{r,q_1}(\mathbb{T}^d)}.$$

**Lemma 24.** **(Fractional Leibniz [12])** *Let $1 < p < \infty$ and $p_1, p_2, p_1', p_2'$ such that*

$$\frac{1}{p_1} + \frac{1}{p_2} = \frac{1}{p_1'} + \frac{1}{p_2'} = \frac{1}{p}.$$

*Then for any $s, \alpha \geqslant 0$ there exists a constant s.t.*

$$\| \langle \nabla \rangle^s (fg) \|_{L^p} \leqslant C \| \langle \nabla \rangle^{s + \alpha} f \|_{L^{p_2}} \| \langle \nabla \rangle^{-\alpha} g \|_{L^{p_1}} + C \| \langle \nabla \rangle^{-\alpha} f \|_{L^{p_2'}} \| \langle \nabla \rangle^{s + \alpha} g \|_{L^{p_1'}}.$$